\documentclass[reqno]{amsart}

\usepackage{graphicx}
\usepackage{hyperref}
\usepackage{ellipsis}
\usepackage{amsmath,amsthm, amssymb}
\usepackage{enumerate}
\usepackage{comment}
\usepackage[font=small,labelfont=bf]{caption}
\usepackage{subcaption}
\usepackage{url}
\usepackage{cite}
\usepackage{color}
\usepackage{color}
\usepackage{breakurl}
\usepackage{comment}
\newcommand{\bburl}[1]{\textcolor{blue}{\url{#1}}}

\numberwithin{equation}{section}

\newtheorem{thm}{Theorem}[section]

\newtheorem{defi}[thm]{Definition}

\theoremstyle{plain}

\newtheorem{remark}[thm]{Remark}

\newcommand\be{\begin{equation}}
\newcommand\ee{\end{equation}}
\newcommand\bea{\begin{eqnarray}}
\newcommand\eea{\end{eqnarray}}
\newcommand\bi{\begin{itemize}}
\newcommand\ei{\end{itemize}}
\newcommand\ben{\begin{enumerate}}
\newcommand\een{\end{enumerate}}
\newcommand\bc{\begin{center}}
\newcommand\ec{\end{center}}
\newcommand\ba{\begin{array}}
\newcommand\ea{\end{array}}

\newcommand{\R}{\ensuremath{\mathbb{R}}}

\newcommand\frakfamily{\usefont{U}{yfrak}{m}{n}}
\DeclareTextFontCommand{\textfrak}{\frakfamily}

\newcommand{\hr}[1]{\href{#1}{\url{#1}}}

\newcommand{\D}{\mathcal{D}}

\title{Crescent configurations}

\author[Burt]{David Burt}
\email{\textcolor{blue}{\href{mailto:drb3@williams.edu}{drb3@williams.edu} }}
\address{Department of Mathematics \& Statistics, Williams College, Williamstown, MA 01267}

\author[Goldstein]{Eli Goldstein}
\email{\textcolor{blue}{\href{mailto:esg2@williams.edu}{esg2@williams.edu} }}
\address{Department of Mathematics \& Statistics, Williams College, Williamstown, MA 01267}

\author[Manski]{Sarah Manski}
\email{\textcolor{blue}{\href{mailto:Sarah.Manski12@kzoo.edu}{Sarah.Manski12@kzoo.edu}}}
\address{Department of Mathematics, Kalamazoo College,  Kalamazoo, MI 49006}

\author[Miller]{Steven J. Miller}
\email{\textcolor{blue}{\href{mailto:sjm1@williams.edu, Steven.Miller.MC.96@aya.yale.edu}{sjm1@williams.edu,Steven.Miller.MC.96@aya.yale.edu} }}
\address{Department of Mathematics \& Statistics, Williams College, Williamstown, MA 01267}

\author[Palsson]{Eyvindur Ari Palsson}
\email{\textcolor{blue}{\href{mailto:eap2@williams.edu}{eap2@williams.edu}}}
\address{Department of Mathematics \& Statistics, Williams College, Williamstown, MA 01267}

\author[Suh]{Hong Suh}
\email{\textcolor{blue}{\href{mailto:hs002012@mymail.pomona.edu}{hs002012@mymail.pomona.edu}}}
\address{Department of Mathematics, Pomona College, Claremont, CA 91711}

\thanks{This work was partially supported by NSF Grants DMS1265673, DMS1347804,  and Williams College. We thank Adam Wang for his help with the Williams College HPCC, and Bill Lenhart and David Moon for many fun conversations about the problem.}



\date{\today}

\begin{document}

\maketitle

\begin{abstract} In 1989, Erd\H{o}s conjectured that for a sufficiently large $n$ it is impossible to place $n$ points in general position in a plane such that for every $1\le i \le n-1$ there is a distance that occurs exactly $i$ times. For small $n$ this is possible  and in his paper he provided constructions for $n\leq 8$. The one for $n=5$ was due to Pomerance while Pal\'{a}sti came up with the constructions for $n=7,8$. Constructions for $n=9$ and above remain undiscovered, and little headway has been made toward a proof that for sufficiently large $n$ no configuration exists. In this paper we consider a natural generalization to higher dimensions and provide a construction which shows that for any given $n$ there exists a sufficiently large dimension $d$ such that there is a configuration in $d$-dimensional space meeting Erd\H{o}s' criteria.
\end{abstract}


\section{Introduction}

In 1946, Erd\H{o}s \cite{Erd46} initiated the first of many problems about distinct distances: what is the minimum number of distinct distances determined by $n$ points in the plane? If we were to randomly place $n$ points in the plane, we would expect all distances between pairs of points to be different. However, as more structure is introduced in the placement of points, some distances may repeat. Erd\H{o}s conjectured that the minimum number of distances is $\Omega(n/\sqrt{\log n})$, which is attained by the $\sqrt{n} \times \sqrt{n}$ integer lattice. The lower bound has been incrementally improved from Erd\H{o}s's original $\Omega(n^{1/2})$ to Larry Guth and Nets Katz's $\Omega(n/\log n)$ announced in 2010, which solves the distinct distances problem up to a $\sqrt{\log n}$ factor \cite{GK}.

Since the introduction of the distinct distances problem, many variants have been tackled. A survey of distance-related problems may be found in \cite{CFG,She}. The variant we study in this paper is that of trying to specify the multiplicity with which each distance occurs. In particular we want each distance to appear a different number of times, with each multiplicity $1$ through $n-1$ represented. This condition affects all possible distances among the $n$ points since there are ${n \choose 2}$ pairs of points and our condition affects $1+2+\ldots+(n-1) = n(n-1)/2$ distances. On a line this problem is trivial, since such a configuration can be achieved by simply placing all $n$ points in an arithmetic progression. In the plane this problem already becomes interesting, provided we prevent too many points from sitting on a single line. Analogously in higher dimensions we want to exclude lower dimensional constructions which we achieve by insisting that the points lie in \emph{general position}.

\begin{defi}[General Position]
We say that $n$ points are in \emph{general position} in $\R^d$ if no $d+1$ points lie on the same hyperplane and no $d+2$ lie on the same hypersphere.
\end{defi}

The multiplicities of the distances are in an increasing order so we introduce the name \emph{crescent configuration} to describe the point configurations we seek.

\begin{defi}[Crescent Configuration]
We say $n$ points are in \emph{crescent configuration (in $\R^d$)} if they lie in general position in $\R^d$ and determine $n-1$ distinct distances, such that for every $1\le i \le n-1$ there is a distance that occurs exactly $i$ times.
\end{defi}

\begin{figure}[h]
\centerline{\includegraphics[height=16em]{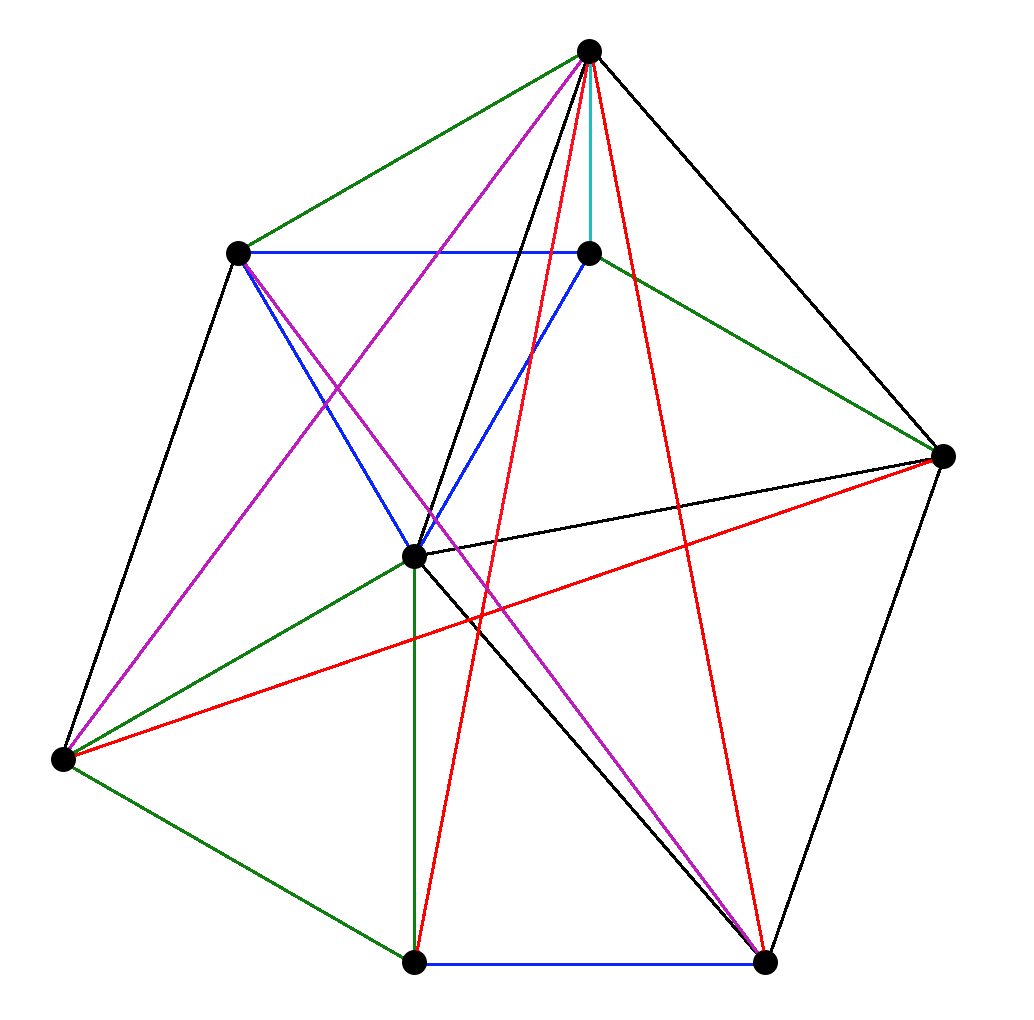}}
\caption{\label{fig:palniseight} Palasti came up with a crescent configuration for $n=8$ \cite{Pal89}. Points have cartesian coordinates $(0,1)$, $(\sqrt{3},0)$, $(2\sqrt{3},0)$, $(\frac{5\sqrt{3}}{2},\frac{5}{2})$, $(\frac{3\sqrt{3}}{2},\frac{9}{2})$, $(\frac{\sqrt{3}}{2},\frac{7}{2})$, $(\frac{3\sqrt{3}}{2},\frac{7}{2})$, $(\sqrt{3},2)$.}
\end{figure}

A one-dimensional variant of this problem of finding crescent configurations, known as the ``beltway" problem, in which all of the points lie on a wrapped interval has music theoretic connections.  In particular, ``deep scales" including commonly used diatonic scales include each interval class a unique number of times\cite{D,Tou}.

In 1989 Erd\H{o}s \cite{Erd89} conjectured that for $n$ sufficiently large, it is not possible to place $n$ points in a crescent configuration in the plane.  Little progress has been made since on this conjecture. Pal\'{a}sti \cite{Pal87,Pal89} provided constructions for $n=7$ and $n=8$, both lying in a small portion of the triangular lattice (see Figure \ref{fig:palniseight}), but constructions for $n=9$ and higher remain undiscovered, with no compelling heuristics to suggest a sufficiently large $n$ such that no configuration exists.

We explore a higher dimensional analogue of this problem. Our main result is the following.

\begin{thm}\label{thm:crescent} For all $n\geq 3$, there exists a set of $n$ points in a crescent configuration in $\R^{n-2}$. \end{thm}

In particular, this shows that given an $n$ there exists a $d$ sufficiently large such that it is possible to place $n$ points in a crescent configuration in $\R^d$.

We give a construction proving Theorem \ref{thm:crescent} in \S\ref{sec:crescentconstr}, and then conclude with some problems for future work in \S\ref{sec:futurework}.


\section{Construction of a crescent configuration}\label{sec:crescentconstr}

We first show by induction that $n-1$ points can be arranged in valid crescent configuration in $\mathbb{R}^{n-2}$.  The base case is clear as an isosceles triangle which is not an equilateral triangle has three points in valid crescent configuration in $\mathbb{R}^2$.

For the inductive step, assume we have $n-2$ points forming a crescent configuration in $\mathbb{R}^{n-3}$.  Any $n-2$ points lie on an $n-4$ dimensional sphere. We now embed the set in $\mathbb{R}^{n-2}$ such that all the distances between points are preserved. There exists a line in $\mathbb{R}^{n-2}$ containing the center of the $n-4$ sphere and perpendicular to the $n-3$ hyperplane containing the sphere.  We can then place the $(n-1)$\textsuperscript{st} point at any distance on this line not previously determined. This new point is equidistant from all of the old points, so it determines a new distance $n-2$ times. Since the other points are in a valid crescent configuration by the inductive hypothesis, they determine distances with the necessary multiplicities from $1$ to $n-3$  Thus, there is an $(n-1)$-point crescent configuration in $\mathbb{R}^{n-2}$, which concludes the inductive argument.

We complete the proof by adding an additional point to an $n-1$ point configuration in $\mathbb{R}^{n-2}$.  This point is added by including the center of the $n-3$ hypersphere defined by the first $n-1$ points.  While in general it is possible that the center of the sphere defines a distance already in the point set, we can easily fix this; we have total control over the radius of the hypersphere through the placement of the $(n-1)$st point on the perpendicular line. Thus, this is a valid $n$-point crescent configuration in $\mathbb{R}^{n-2}$ and completes the proof. \hfill $\Box$

\section{Future Work}\label{sec:futurework}

It seems likely that the actual dimension needed to place points in a valid configuration is significantly lower. There are a number of other related problems of interest. Define $\D(n)$ to be the minimum dimension, greater than 1, such that $n$ points can be placed in $\R^{\D(n)}$ in a crescent configuration. Our construction shows that $\D(n)\leq n-2$ for all $n>3$.\\ \

\begin{itemize}
\item Albujer asked, is $\D(n)$ bounded as $n$ goes to infinity \cite{Alb}?\\ \

\item Is $\D(n)$ sublinear? \\ \

\item Is $\D(n)$ monotonically increasing? \\ \

\item If $\D(n)=k$ is it possible to place $n$ points in a valid configuration in $\R^d$ for all $d>k$? Note that in general we cannot just embed the solution into higher dimensions and maintain general position. For example, if you take the 8 point construction from $\mathbb{R}^2$ by Pal\'{a}sti and embed it in the most natural way into $\mathbb{R}^3$ then immediately you have more than $4$ points in a plane which violates the condition of general position. \\ \ 

\item Do any planar constructions for $n \geq 9$ exist? This has been asked before, see for example \cite{CFG,Alb}.\\ \ 

\item Can planar constructions for $n\geq 9$ be found on the triangular lattice?  It is known that constructions for $n<9$ exist on the triangular lattice. \\ \

\end{itemize}

\begin{remark}
With the help of a parallel computing cluster, we have exhaustively searched a 91 point hexagonal region of the triangular lattice for a construction for $n=9$, but none exist. As the naive implementation took over 900 hours of computation for this size, better (and achievable) techniques are required to search a substantively larger region.
\end{remark}


\ \\
\end{document}